\documentclass[11pt,leqno]{amsart}
\usepackage{amsmath,amsfonts,amssymb,amscd,amsthm,amsbsy,upref}
\newtheorem{thm}{Theorem}[section]
\newtheorem{lem}[thm]{Lemma}
\newtheorem{cor}[thm]{Corollary}
\newtheorem{prop}[thm]{Proposition}
\theoremstyle{definition}
\newtheorem{defn}[thm]{Definition}
\newtheorem{prob}[thm]{Problem}    
\newtheorem*{rem}{Remark}
\newtheorem{remark}[thm]{Remark}
\def\nat{{\mathbb N}}
\def\real{{\mathbb R}}
\def\A{{\mathcal A}}
\def\F{{\mathcal F}}

\def\bD{{\mathbf D}}
\def\bE{{\mathbf E}}
\def\bF{{\mathbf F}}
\def\bH{{\mathbf H}}
\def\bG{{\mathbf G}}
\def\ep{\varepsilon}
\def\w{\omega}
\def\FDD{\operatorname{FDD}}
\def\trivert#1{|\!|\!| #1|\!|\!|}
\def\Bigtrivert#1{\Big|\!\Big|\!\Big| #1\Big|\!\Big|\!\Big|}
\def\th{\tilde h}
\def\tg{\tilde g}
\def\tw{\tilde w}
\def\tz{\tilde z}
\def\tE{\widetilde E}
\def\tbE{\widetilde{\bE}}
\def\tH{\widetilde H}
\def\tbH{\widetilde{\bH}}
\def\tQ{\widetilde Q}
\def\tZ{\widetilde Z}
\def\oK{\overline{K}}

\begin{document}
\title{A universal reflexive space for the class of uniformly convex 
Banach spaces}
\author{E. Odell \and Th. Schlumprecht}
\address{Department of Mathematics\\ The University of Texas at Austin\\
1 University Station C1200\\ Austin, TX 78712-0257}
\email{odell@math.utexas.edu}
\address{Department of Mathematics\\ Texas A\&M University\\
College Station, TX 77843-3368
}
\email{schlump@math.tamu.edu}
\dedicatory{Dedicated to the memory of V. I. Gurarii}
\thanks{Research supported by the National Science Foundation.}

\begin{abstract}
We show that there exists a separable reflexive Banach space  into
which every separable uniformly convex Banach space isomorphically embeds. 
This solves  problems raised by  J.~Bourgain \cite{B} in 1980 and
 by W.~B.~Johnson in 1977 \cite{Jo}. 
We also give intrinsic characterizations of separable reflexive Banach 
spaces which embed into a reflexive space with a block $q$-Hilbertian 
and/or a block $p$-Besselian finite dimensional decomposition.
\end{abstract}

\maketitle

\section{Introduction}

J.~Bourgain \cite{B} proved that if $X$ is a separable Banach space which 
contains an isomorph of every separable reflexive space then $X$ contains 
an isomorph of $C[0,1]$ and hence is universal, i.e., $X$ contains an 
isomorph of every separable Banach space. 
He asked if there exists a separable reflexive space $X$ which is universal 
for the class of all separable uniformly convex (equivalently, all 
superreflexive \cite{E}, \cite{Pi}) Banach spaces. 
Such an $X$ could not be superreflexive since $c_0$ and $\ell_1$ are finitely 
representable in any space which contains isomorphs of all $\ell_p$'s for 
$1<p<\infty$. 

We shall answer Bourgain's question in the affirmative. 
S.~Prus \cite{P} gave a partial solution by proving that there exists a 
reflexive Banach space $X$ which is universal for all spaces with a finite 
dimensional decomposition (FDD) which satisfy  $(p,q)$-estimates for some 
$1<q\le p<\infty$. 

\begin{defn}	
Let $(F_n)$ be an FDD. 
$(x_n)$ is a {\em block sequence} of $(F_n)$ if there exist integers 
$0=k_0 < k_1 <\cdots$ so that for all $n\in\nat$,
$$x_n \in [F_i]_{i\in (k_{n-1},k_n]} 
\equiv \text{span} \{F_i : k_{n-1} <i \le k_n\}\ .$$
\end{defn}

\begin{defn}	
Let $1\le q\le p\le\infty$ and let $C<\infty$. 
An FDD $(F_n)$ satisfies {\em $C$-$(p,q)$-estimates} if for all block 
sequences $(x_n)$ of $(F_n)$, 
$$C^{-1} \left(\sum \|x_n\|^p\right)^{1/p} \le \|\sum x_n\| 
\le C\left(\sum \|x_n\|^q\right)^{1/q}\ .$$
\end{defn}

We say that {\em $(F_n)$ satisfies $(p,q)$-estimates} if it satisfies 
$C$-$(p,q)$-estimates for some $C<\infty$. 
A {\em basic sequence $(x_n)$ is  said to satisfy $(p,q)$-estimates}
if $(E_n)$ does where 
$E_n = \text{span} \{x_n\}$ for $n\in\nat$.

\subsection*{Terminology}
In some of the literature an FDD satisfying $(p,1)$-estimates is called 
block $p$-Besselian and one satisfying $(\infty,q)$-estimates is called 
block $q$-Hilbertian.

We shall prove that if $X$ is uniformly convex then there exists 
$1<q\le p<\infty$ and a space $Z$ with an FDD satisfying $(p,q)$-estimates 
such that $X$ embeds into $Z$. 
In combination with Prus' result we then obtain the solution to 
Bourgain's problem.

\begin{thm}\label{Bourgain-solution}
There exists a separable reflexive Banach space $X$ 
which contains an isomorph of 
every separable superreflexive Banach space.
\end{thm}

To accomplish this we shall characterize when a reflexive space embeds into 
a reflexive space with an FDD satisfying $(p,q)$-estimates. 
Before stating our results in this regard we need some more definitions.

\begin{defn} 	
If $\bE = (E_n)$ is an FDD for a space $X$, by $P_n^\bE$ we denote 
the natural projection of $X$ onto $E_n$. 
More generally if $I$ is an interval or finite union of intervals in $\nat$, 
$P_I^\bE$ shall denote the natural projection on $X$ given by 
$P_I^\bE (\sum e_n) = \sum_{n\in I} e_n$ (where $e_n\in E_n$ for all $n$). 
The {\em projection constant} of $(E_n)$ is $\sup \{\| P_I^\bE\| : I$ is an 
interval in $\nat\}$. 
$(E_n)$ is {\em bimonotone} if its projection constant is~1. 
A {\em blocking} $(G_n)$ of $(E_n)$ is an FDD given by 
$G_n = [E_i]_{i\in (N_{n-1},N_n]}$ for some sequence of integers 
$0 = N_0 < N_1 < N_2 < \cdots$.
\end{defn}

Henceforth all Banach spaces will be assumed to be separable. 
$S_X$ denotes the unit sphere of $X$ and $B_X$ denotes the unit ball of $X$.

\begin{defn}	
a) $T_\infty \equiv \{(n_1,\ldots,n_k) : k\in \nat$ and $n_1 < n_2 < \cdots 
< n_k$ are natural numbers$\}$. 
$T_\infty$ is  ordered by $(n_1,\ldots,n_k) \le 
(m_1,\ldots,m_\ell)$ iff $k\le \ell$ and $n_i = m_i$ for $i\le k$.

b)  A {\em tree in } a Banach space $X$ is a family in $X$ indexed by 
 $T_\infty$. 
A {\em weakly null tree} in $X$ is a tree 
  $(x_\alpha)_{\alpha\in T_\infty} \subseteq X$ with the property that for all 
$\alpha=(n_1,n_2,\ldots n_k)\in T_\infty\cup \{\emptyset\}$, 
  $(x_{(\alpha,n)})_{n>n_k}^\infty$ is weakly null. 
$(y_i)_{i=1}^\infty$ is a {\em branch} of $(x_\alpha)_{\alpha\in T_\infty}$ 
if there exist $n_1 < n_2 < \cdots$ so that $y_k = x_{(n_1,\ldots,n_k)}$ 
for all $k\in\nat$.

c) If $(x_\alpha)_{\alpha\in T_{\infty}}$ is a tree and if $T'\subset T_\infty$
is such that 	
for each $\alpha\in T'\cup \{\emptyset\}$ there
is an infinite $N_\alpha\subset\nat$ so that 
$(\alpha,n)\in T'$ for all $n\in N_\alpha$ we call 
$(x_\alpha)_{\alpha\in T'}$ a {\em full subtree}. 
In this case we can relabel   
$(x_\alpha)_{\alpha\in T'}$ into $(z_\alpha)_{\alpha\in T_\infty}$.
Note that 	
every branch of a full subtree is a branch of the original tree. 

\end{defn}

\begin{defn}	
Let $1\le q\le p\le\infty$ and $C<\infty$. 
A Banach space {\em $X$ satisfies $C$-$(p,q)$-tree estimates} if for all 
weakly null trees in $S_X$ there exists  branches $(y_i)$ and $(z_i)$ 
satisfying 
$$C^{-1} \left( \sum |a_i|^p\right)^{1/p} 
\le \|\sum a_i y_i\| 
\ \text{ and }\ \|\sum a_i z_i\| 
\le C\left( \sum |a_i|^q\right)^{1/q}$$
for all $(a_i) \subseteq \real$. 
{\em $X$ satisfies $(p,q)$-tree estimates} if it satisfies 
$C$-$(p,q)$-tree estimates for some $C<\infty$.
\end{defn}

\begin{thm}\label{mainthm}
Let $X$ be a reflexive Banach space and let $1\le q\le p\le\infty$. 
The following are equivalent.
\begin{itemize}
\item[{a)}] $X$ satisfies $(p,q)$-tree estimates.
\item[{b)}] $X$ is isomorphic to a subspace of a reflexive space $Z$ having 
an $\FDD$ which satisfies $(p,q)$-estimates.
\item[{c)}] $X$ is isomorphic to a quotient of a reflexive space $Z$ having 
an $\FDD$ which satisfies $(p,q)$-estimates.
\end{itemize}
\end{thm}

Theorem~\ref{Bourgain-solution}  is a corollary of this 
(using Prus' result \cite{P}) since for every uniformly convex space $X$ 
there exists $K<\infty$ and $1<q\le p<\infty$ such that every normalized 
2-basic sequence in $X$ admits $K$-$(p,q)$-estimates (\cite{J}, \cite{GG}). 
Indeed it is trivial to extract a 2-basic branch from a normalized weakly null tree. 
Theorem~\ref{mainthm} also solves problem~IV.3 in \cite{Jo}. 

\section{The Proof}

The equivalence of a) and b) in Theorem~\ref{mainthm} in the case 
$1< q = p<\infty$ was established in \cite{OS}. 
We shall be using some blocking arguments established there and in earlier 
seminal
papers of W.B.~Johnson and M.~Zippin (\cite{Jo}, \cite{JZ1,JZ2}) which we 
shall recall as needed. 
A key first step will of course be Zippin's result \cite{Z} that a reflexive 
space embeds into a reflexive space with an FDD (in fact with a  basis). 

Before stating Theorem~\ref{main-central}, which contains the central part 
of our main Theorem~\ref{mainthm}, we set some more notation. 
If $(E_n)$  is an FDD then by $c_{00} (\oplus_{n=1}^\infty E_n)$ we mean 
the subspace of all $x= \sum e_n$ where $e_n \in E_n$ for all $n$ and only 
finitely many $e_n$'s are nonzero. 
If $Z$ has an FDD, $\bF = (F_n)$, and $1<p<\infty$ then $Z_p (\bF)$ denotes 
the Banach space obtained by completing $c_{00} (\oplus_{n=1}^\infty F_n)$ 
under $\|\cdot \|_{_{Z_p}}$ given by: for $y = \sum y_n$, $y_n \in F_n$ for all $n$, 
$$\|y\|_{_{Z_p}} = \sup \left\{ \bigg( \sum_{j=1}^\infty 
\Big\| \sum_{i=n_{j-1} +1}^{n_j} y_i\Big\|^p\bigg)^{1/p} 
: 0 = n_0 < n_1 < \cdots\right\}\ .$$
Note that $(F_n)$ is a bimonotone FDD for $Z_p (\bF)$ satisfying 
1-$(p,1)$-estimates. 

\begin{thm}\label{main-central}
Let $X$ be a reflexive Banach space and let $1<p<\infty$. 
If $X$ satisfies $(p,1)$-tree estimates then 
\begin{itemize}
\item[{a)}] $X$ can be embedded into a reflexive space $Z$ with an $\FDD$ 
satisfying $(p,1)$-estimates.
\item[{}] $\quad$ More precisely, if $X$ is a subspace of 
$Z$, a reflexive space with an 
$\FDD\ (E_n)$ then there exists a blocking $\bF= (F_n)$ of $(E_n)$ so that $X$ 
naturally embeds into the reflexive space $Z_p(\bF)$.
\item[{b)}] $X$ is the quotient of a reflexive space with an $\FDD$ satisfying 
$(p,1)$-estimates.
\end{itemize}
\end{thm}

The proof of a) is much like the proof in \cite{OS}. 
The proof of b) requires some new ideas. 
Before starting the proof we need some terminology and preliminary results.

\begin{defn}	
Let $\bE = (E_i)$ be an FDD for $Y$ and let $\delta =  (\delta_i)$ with 
$\delta_i\downarrow 0$. 
A sequence $(y_i) \subseteq S_Y$ is called a {\em $\delta$-skipped block w.r.t.
$(E_n)$} if there exist integers $1= k_0 < k_1 < \cdots$ so that for all 
$i\in \nat$, 
$$\|P_{(k_{i-1},k_i)}^\bE y_i - y_i\| < \delta_i\ .$$
\end{defn}

\begin{defn}	
If $\A \subseteq S_X^\w$, the set of all normalized sequences in $X$, and 
$\ep >0$ we set  
$$\A_\ep = \{(x_n)\in S_X^\w :\text{ there exists } (y_n) \in \A
\text{ with } \|x_n-y_n\| < \frac{\ep}{2^n}\text{ for all } n\}\ .$$
\end{defn}

\noindent
$\overline{\A_\ep}$ denotes the closure 
of $\A_\ep$ w.r.t.\ the product topology of the 
discrete topology on $S_X$. 

The next result is Theorem~3.3 b)~$\Leftrightarrow$ d) in 
\cite{OS}. 

\begin{prop}\label{OSthm3.3}
Let $X$ be a Banach space with a separable dual. Then  $X$ is (isometrically) a 
subspace of a 
  Banach space $Z$ having a shrinking FDD $(E_n)$ satisfying the 
following:

 For $\A\subseteq S_X^\w$.  the following are equivalent.
\begin{itemize}
\item[{a)}] For all $\ep>0$ every weakly null tree in $S_X$ has a branch in 
$\overline{\A_\ep}$. 
\item[{b)}]  For all $\ep>0$ there exists a blocking $(F_i)$ of $(E_i)$ 
and $\delta = (\delta_i)$, $\delta_i\downarrow 0$, so that if 
$(x_n) \subseteq S_X$ is a $\delta$-skipped block w.r.t.\ $(F_i)$ then 
$(x_n) \in \overline{\A_\ep}$.
\end{itemize}
\end{prop}

The following Proposition yields that
in the reflexive case the equivalence (a)$\iff$(b) in Proposition
 \ref{OSthm3.3}  holds for any embedding of $X$ into a reflexive Banach space 
$Z$ with an FDD.

\begin{prop}\label{addOSthm3.3} Let $Z$ and $Y$ be reflexive  spaces with FDDs 
$\bE=(E_n)$ and $\bF=(F_n)$,
 respectively, both containing a space $X$,  and let 
$\delta=(\delta_n)\subset (0,1)$, with
 $\delta_n\downarrow 0$, as $n\uparrow \infty$.

Let $C$ be the maximum of the projection  constants of $(E_n)$ and 
$(F_n)$.
 Then there is a blocking $\bG=(G_n)$ of $(F_n)$, so that every normalized 
$\frac{\delta}{5C^3}$-skipped
 block of $(G_n)$ in $S_X$  is   a $\delta$-skipped block of $(E_n)$.
\end{prop}  

\begin{proof} By induction we will choose  $0=M_0<M_1<M_2<\ldots $
 and $N_1<N_2<\ldots$ in $\nat$ so that for all $k\in\nat$
 \begin{align}
\label{E:2.4.2}
&\forall x\in S_X \,\forall i\in\{1,2\ldots ,k,k+1\}\,\,
\text{ if }\|P^{\bF}_{(M_{k-1},\infty)}(x)\|\le \frac{\delta_i}{5C^2}\text{ 
then }
           \|P^{\bE}_{[N_k,\infty)}(x)\|\le \frac{\delta_i}{2},\\
\label{E:2.4.1}
&\forall x\in S_X \,\forall i\in\{1,2\ldots k,k+1\}\,\, 
\text{ if }\|P^{\bF}_{[1,M_k]}(x)\|\le \frac{\delta_i}{5C^2}\text{ then }
           \|P^{\bE}_{[1,N_k]}(x)\|\le \frac{\delta_i}{2}
\end{align}
Once  accomplished  we choose $G_k=\oplus_{i=M_{k-1}+1}^{M_k} F_i$.
 If $(x_n)$ is a $\delta/ {5C^3}$-skipped
 block of $(G_n)$ in $S_X$,  there exist $0=k_0<k_1<k_2<\ldots $  such that 
for all
 $n\in\nat$
$$\|x_n-P^{\bG}_{(k_{n-1},k_n)}(x_n)\|
 =\|x_n-P^{\bF}_{(M_{k_{n-1}},M_{k_n})}(x_n)\|\le\frac{\delta_n}{5C^3}.$$
Thus,
$$\| P^{\bF}_{[1,M_{k_{n-1}}]}(x_n)\|\le\frac{\delta_n}{5C^2}\text{ and }
 \| P^{\bF}_{(M_{k_n-1},\infty)}(x_n)\| \le\frac{\delta_n}{5C^2}.$$
We  deduce from (\ref{E:2.4.1})  and (\ref{E:2.4.2}) that
$$ \| P^{\bE}_{[1,N_{k_{n-1}}]} (x_n)\|\le \frac{\delta_n}{2}\text{ and }
 \| P^{\bE}_{[N_{k_n},\infty)}(x_n)\|\le \frac{\delta_n}{2},$$
which yields that $(x_n)$ is a $\delta$-skipped block of $(E_i)$.

Assume that we have chosen $M_{k-1}$ for some $k\ge 1$.  
 We need to find an $N_k$ which satisfies
 (\ref{E:2.4.2}). If such an $N_k$ did not exist,  we could find sequences 
$(x_j)\subset S_X$  and
$(i_j)\subset\{1,2\ldots k+1\}$ so that for any $j > N_{k-1}$, 
$$\|P^{\bF}_{(M_{k-1},\infty)}(x_j)\|\le \frac{\delta_{i_j}}{5C^2}\text{ and  }
          \|P^{\bE}_{[j,\infty )}(x_j)\|>\frac{\delta_{i_j}}{2}.$$
Passing to a subsequence  we may assume that $i_j=i$ for all 
$j\in J$ and some 
$i\in\{1,2\ldots k+1\}$, where $J$ is a subsequence of $\nat$. 
Since $\lim_{j\to\infty,\, j\in J} \|P^{\bE}_{[j,\infty)}\circ 
P^{\bF}_{(M_{k-1},\infty)}-P^{\bE}_{[j,\infty)}\|=0$
we deduce that
 \begin{align*}
\frac{\delta_i}2&\le \limsup_{\substack{j\to\infty \\ j\in J}}
\|P^{\bE}_{[j,\infty)}(x_j)\|\le
  \limsup_{\substack{j\to\infty\\ j\in J}}
  \|P^{\bE}_{[j,\infty)}\circ 
P^{\bF}_{(M_{k-1},\infty)}(x_j)\|\\
 &\le C
\limsup_{\substack{j\to\infty\\ j\in J}}\|P^{\bF}_{(M_{k-1},\infty)}(x_j)\|\le 
\frac{\delta_{i}}{5C^2},
\end{align*}
which is a contradiction, and finishes the proof of our claim.

Assume now that we have chosen $N_k$,  but there is no $M_k$ satisfying
(\ref{E:2.4.1}).  We could choose a sequence $(x_j)\subset S_X$  and
$(i_j)\subset\{1,2\ldots k+1\}$ so that for any $j> M_{k-1}$
$$\|P^{\bF}_{[1,j]}(x_j)\|\le \frac{\delta_{i_j}}{5C^2}\text{ and }
\|P^{\bE}_{[1,N_k]}(x_j)\|> \frac{\delta_{i_j}}{2}.$$
After passing to subsequences we can assume that $i_j=i$ for some fixed 
$i\in\{1,2,\ldots k+1\}$ and $j\in J$, a subsequence of $\nat$, 
and that $(x_j)_{j\in J}$ converges weakly to some $x\in B_X$.
Then it follows that
\begin{align*}
\|x\|&=\lim_{j_0\to\infty}\| P^{\bF}_{[1,j_0]}(x)\| 
=\lim_{j_0\to\infty} \lim_{\substack{j\to\infty\\ j\in J}}  
\| P^{\bF}_{[1,j_0]}(x_j)\|\\
&\le \limsup_{j_0\to\infty} \limsup_{\substack{j\to\infty\\ j\in J}} 
\| P^{\bF}_{[1,j]}(x_j)\| + \|P^{\bF}_{(j_0,j]}(x_j)\|\\
&\le (1+C) \limsup_{\substack{j\to\infty\\ j\in J}}\| P^{\bF}_{[1,j]}(x_j)\| 
\le \frac{2\delta_i}{5C}
\end{align*}
and that
$$\|x\|\ge \frac1C\|P^{\bE}_{[1,N_k]}(x)\|
=\frac1C\lim_{\substack{j\to\infty\\ j\in J}} 
\|P^{\bE}_{[1,N_k]}(x_j)\|\ge 
\frac{\delta_i}{2C},$$
which is a contradiction.
\end{proof}

{From} Corollary 4.4 in \cite{OS} we have 

\begin{prop}\label{OSprop4.4}
Let $X$ be a Banach space which is a subspace of a reflexive space $Z$ with 
an $\FDD\ \bE = (E_i)$ having projection constant $K$. 
Let $\delta_i\downarrow 0$. 
Then there is a blocking $\bF= (F_i)$ of $(E_i)$ given by 
$F_n = [E_i]_{i\in (N_{n-1},N_n]}$ 
for some integers $0= N_0 < N_1 < \cdots$ with the following property. 
For all $x\in S_X$ there exists $(x_i) \subseteq X$ and integers $(t_i)$ 
with $t_i\in (N_{i-1},N_i)$ for all $i$ such that 
\begin{itemize}
\item[{a)}] $x = \sum\limits_{i=1}^\infty x_i$
\item[{b)}] For $i\in\nat$ either $\|x_i\| <\delta_i$ or 
$\|P_{(t_{i-1},t_i)}^\bE x_i - x_i\| < \delta_i \|x_i\|$.
\item[{c)}] $\|P_{(t_{i-1},t_i)}^\bE x-x_i\| < \delta_i$ for all $i\in\nat$.
\item[{d)}] $\|x_i\| < K+1$ for $i\in\nat$.
\item[{e)}] $\|P_{t_i}^\bE x\| < \delta_i$ for $i\in\nat$.
\end{itemize}
Moreover the above hold for any further blocking of $(F_n)$ (which would 
redefine the $N_i$'s). 
\end{prop}

Parts d) and e) were not explicitly stated in \cite{OS} but follow from 
the proof.

\begin{proof}[Proof of Theorem~\ref{main-central} a)]
Let $X$ be contained in a reflexive space $Z$ with an $\FDD\ \bE=(E_i)$ having 
projection constant $K$. 
Assume that $X$ satisfies $C$-$(p,1)$-tree estimates. 
Let $\A = \{(x_i) \in S_X^\w$:  $(x_i)$ is $\frac32$-basic and for all scalars 
$(a_i)$, $C\|\sum a_i x_i\| \ge ( \sum |a_i|^p)^{1/p}\}$.
Choose $\ep>0$ so that if $(x_i) \in \overline{\A_\ep}$ then $(x_i)$ 
is 2-basic and satisfies for all $(a_i)\subseteq \real$, 
$$2C\|\sum a_i x_i\| \ge \left( \sum |a_i|^p\right)^{1/p}\ .$$

By Propositions~\ref{OSthm3.3} and \ref{addOSthm3.3} there exists $\delta = (\delta_i)$, 
$\delta_i \downarrow 0$, and a blocking of $(E_i)$, which we still denote 
by $(E_i)$, so that every $\delta$-skipped block w.r.t.\ $(E_i)$ is in 
$\overline{\A}_\ep$. 
We then use $\bar \delta = (\bar\delta_i)$ where 
$\bar\delta_i = \delta_i/2K$ to form a new blocking 
$F_n = [E_i]_{i\in (N_{n-1},N_n]}$ 
satisfying the conclusion of Proposition~\ref{OSprop4.4}.
We assume, as we may, that $\sum_{i=1}^\infty \delta_i < 1/3$.
Note that any subsequence of a $\bar\delta$-skipped block w.r.t.\ 
$(F_i)$ is then a $\delta$-skipped block w.r.t.\ $(E_i)$.

Our goal is to prove that $X$ naturally embeds into $Z_p (\bF)$. 
To achieve this we prove that if $x\in S_X$ then, for some absolute 
constant $A= A(K,C)$, $(\sum \|P_n^\bF x\|^p)^{1/p} \le A$. 
Since the argument we will give would also work for any blocking of $(F_n)$ 
(see the ``moreover'' part of Proposition~\ref{OSprop4.4}) we obtain 
$\|x\|_{_{Z_p}} \le A$ which finishes the proof of the claim.

Let $x\in S_X$ and write $x= \sum x_i$ with $(x_i) \subseteq X$ and 
$t_i\in (N_{i-1},N_i]$ as in Proposition~\ref{OSprop4.4}. 
Let $y_i = P_{(t_{i-1},t_i]}^\bE x$ for $i\in\nat$.
Let $B = \{i\ge 2: \|x_i\| \ge \bar\delta_i$ and $\|P_{(t_{i-1},t_i)}^\bE 
x_i - x_i\| < \bar\delta_i \|x_i\|\}$. 
Since $(x_i/\|x_i\|)_{i\in B}$ is a $\delta$-skipped block w.r.t.\ $(E_i)$
we know it is in $\overline{\A_\ep}$ and so 
$2C\|\sum_{i\in B} x_i\| \ge (\sum_{i\in B}  \|x_i\|^p)^{1/p}$. 
Also for all $i\in \nat$,
$$\|y_i\| = \|P_{(t_{i-1},t_i)}^\bE x-x_i + P_{t_i}^\bE x + x_i\| 
\le 2\bar\delta_i + \|x_i\|$$
by c) and e) of Proposition~\ref{OSprop4.4}. 
Thus $\|y_1\| < K+2$ and 
\begin{equation*}
\begin{split}
\sum \|y_i\|^p 
& = \sum_{i\in B} \|y_i\|^p + \sum_{i\notin B} \|y_i\|^p 
\le \sum_{i\in B} (2\bar\delta_i + \|x_i\|)^p  + \|y_1\|^p
+ \sum_{i\notin B} (3\bar\delta_i)^p\\
& < \sum_{i\in B} 3^p \|x_i\|^p + (K+2)^p + 1 
\le 3^p (2C)^p \|\sum_{i\in B} x_i\|^p + (K+2)^p +1\ .
\end{split}
\end{equation*}
Now 
$$\|\sum_{i\in B} x_i\| \le 1 + \sum_{i\notin B} \|x_i\| 
< 1 + \|x_1\| + \sum_i \bar \delta_i < K+3 \ .$$
Thus 
$$\sum \|y_i\|^p \le (6C)^p (K+3)^p + (K+2)^p +1\equiv A'\ .$$
Let $z_i = P_{(N_{i-1},N_i]}^\bE x = P_i^\bF x = P_i^\bF (y_i + y_{i+1})$ 
for $i\in\nat$. 
Thus $\|z_i\| \le K (\|y_i\| +  \|y_{i+1}\|)$ and so 
$$\left(\sum \|z_i\|^p\right)^{1/p}
\le 2K [A']^{1/p} \equiv A\ .$$
\renewcommand{\qed}{}
\end{proof}

To complete the proof of part a) we have the following easy 

\begin{lem}\label{part-a)}
Let $\bF= (F_i)$ be a shrinking $\FDD$ for a Banach space $Z$. 
Then for $1<p<\infty$, $Z_p (\bF)$ is reflexive.
\end{lem}

\begin{proof}
As noted earlier, $(F_i)$ is a bimonotone FDD for $Z_p(\bF)$ which satisfies 
1-$(p,1)$-estimates and hence $(F_i)$ is boundedly complete. 
Let  $\frac1p + \frac1{p'} =1$  and set
$$\F = \left\{\sum a_i f_i : \begin{matrix} \text{$(a_i)\in B_{\ell_{p'}}$ and $(f_i)$ is a (finite or} \\
                                            \text{infinite) block sequence of $(F_n^*)$ in $S_{Z^*}$}
                                                \end{matrix}
                                      \right\}.$$ 
 If the above sum $\sum a_i f_i$ is a finite 
one, say $\sum_{i=1}^n a_i f_i$, then $f_n$ can be supported on 
$[F_m^*]_{m\in [j,\infty)}$ for some $j$. 
It is easy to check that $\F$ is a weak* compact 1-norming subset 
of $Z_p (\bF)^*$. 
Thus $Z_p(\bF)$ is isometrically a subspace of $C(\F)$, the space of continuous
 function on $\F$.
 Since every $\|\cdot\|_{_{Z_p}}$-normalized block $(z_i)$ of $(F_i)$ is 
pointwise null on $\F$, hence weakly null, it follows that $(F_i)$ is a shrinking FDD for $Z_p(\bF)$.
\end{proof}

To prove part b) we need a blocking result due to Johnson and Zippin.

\begin{prop}\label{JZresult}
\cite{JZ1}
Let $T:Z\to W$ be a bounded linear operator from a space $Z$ with a 
shrinking $\FDD\ (G_n)$ into a space $W$ with an $\FDD\ (H_n)$. 
Let $\ep_i\downarrow 0$. 
Then there exist blockings $\bE = (E_n)$ of $(G_n)$ and 
$\bF = (F_n)$ of $(H_n)$ so that: 
for all $i\le j$ and $z\in S_{[E_n]_{n\in (i,j]}}$ we have 
$\|P_{[1,i)}^\bF Tz\| < \ep_i$ 
and $\|P_{(j,\infty)}^\bF Tz\| <\ep_j$.
\end{prop}

\begin{proof}[Proof of Theorem~\ref{main-central} b)]
By Lemma~3.1 in \cite{OS} we can, by renorming, regard $X^* \subseteq Z^*$ 
where $Z^*$ is a reflexive space with a bimonotone $\FDD\ (E_i^*)$ such that 
$c_{00} (\oplus_{i=1}^\infty E_i^*) \cap X^*$ is dense in $X^*$. 
Thus we have a quotient map $Q :Z\to X$. 
By part~a) we may regard $X\subseteq W$, a reflexive space with an 
$\FDD\ (F_i)$ satisfying $C$-$(p,1)$-estimates for some $C$. 
Let $K$ be the projection constant of $(F_i)$.

Choose $\delta = (\delta_i)$, $\delta_i\downarrow 0$, so that if $(y_i)$ 
is any $\delta$-skipped block of any blocking of $(F_i)$ and $(z_i)$ satisfies 
$\|z_i -y_i\| < 3K\delta_i$ for all $i$, then $(z_i)$ is 2-equivalent to 
$(y_i)$, is $2K$-basic and $(y_i)$ satisfies $2C$-$(p,1)$-estimates. 

In addition we require that 
\begin{equation}\label{eq1}
\sum_{i=1}^\infty 3(K+1) K\delta_i < \frac14\ ,\text{ and }
6/(1-3K\delta_1) < 7
\end{equation}
and we choose $\ep_i \downarrow 0$ with $6\ep_i <\delta_i$ for all $i$.

By Proposition~\ref{JZresult}, blocking and relabeling our FDD's we may assume
\begin{equation}\label{eq2}
\text{For all $i\le j$ and $z\in S_{[E_n]_{n\in (i,j]}}$ we have}
 \|P_{[1,i)}^\bF Qz\| < \ep_i\ \text{ and }\ 
\|P_{(j,\infty)}^\bF Qz\| <\ep_j\ .
\end{equation}

For $i\in\nat$ let $\widetilde E_i$ be the quotient space of $E_i$ determined 
by $Q$. 
Thus if $z\in E_i$, the norm on $\tz $ (the equivalence class of $z$ in 
$E_i$) is $\trivert{\tz } = \|Qz\|$. 
We may assume $\widetilde E_i\ne \{0\}$ for all $i$. 
More generally for $\tz  = \sum \tz _i \in c_{00} (\oplus_{i=1}^\infty
\widetilde E_i)$ with $\tz _i \in \widetilde E_i$ for every $i$, we set 
$$\trivert{\tz } = \sup_{m\le n} \Big\| \sum_{i=m}^n Qz_i\| 
= \sup_{m\le n} \|QP_{[m,n]}^\bE z\|\ .$$

We let $\tZ $ be the completion of $(c_{00} (\oplus_{i=1}^\infty 
\widetilde E_i),\trivert{\cdot})$. 
Note that if $\tz  = \sum \tz _i \in c_{00} (\oplus_{i=1}^\infty 
\widetilde E_i)$ then setting $\tQ \tz  \equiv \sum \tQ
\tz _i \equiv \sum Qz_i$, we have $\|\tQ \tz \| \le 
\trivert{\tz }$. 
Thus $\tQ$ extends to a norm one map from $\tZ $ into $X$.
Before continuing the proof of b) we need 
\renewcommand{\qed}{}
\end{proof}

\begin{prop}\label{shrinkingFDD}
$\quad$
\begin{itemize}
\item[{a)}] $(\tilde E_i)$ is a bimonotone shrinking $\FDD$ for $\tZ $.
\item[{b)}]  $\tQ$ is a quotient map from $\tZ $ onto $X$. 
More precisely if $x\in X$ and $z\in Z$ with $Qz = x$, $\|z\| = \|x\|$, 
and $z = \sum z_i$
with $z_i\in E_i$, then $\tz  = \sum \tz _i \in \tZ$, 
$\trivert{\tz } = \|z\|$ and $\tQ\tz  = x$.
\item[{c)}] Let $(\tz _i)$ be a block sequence of $(\tilde E_i)$ in 
$B_{\tZ }$ and assume that $(\tQ\tz _i)$ is a basic sequence 
with projection constant $\overline{K}$ and 
$a= \inf_i \|\tQ\tz _i\|>0$.
Then for all scalars $(a_i)$ we have 
\begin{equation*}
\|\sum a_i \tQ(\tz _i) \| 
 \le \trivert{\sum a_i\tz _i} 
 \le \frac{3\overline{K}}{a} \|\sum a_i \tQ \tz _i\|\ .
\end{equation*}
\end{itemize}
\end{prop}

\begin{proof}
By definition $\tilde{\bE} =(\tilde E_i)$ is a bimonotone FDD for $\tZ $. 
We will deduce later that it is shrinking.
To see b) let $Qz = x$ with $\|z\| = \|x\|$ and write $z = \sum z_i$ 
with $z_i \in E_i$ for all $i$. 
Then for $i\le j$, 
$$\trivert{\sum_{\ell=i}^j \tilde z_\ell}  
=\sup_{i\le n\le m\le j} \| Q(\sum_{\ell=n}^m z_\ell)\| 
\le \sup_{i\le n\le m\le j} \|\sum_{\ell=n}^m z_\ell\| 
=\|\sum_{\ell=i}^j z_\ell\|\ .$$
Thus $\sum \tz_i$ converges in $\tZ$ to some $\tz$ with, $\trivert{\tz} = \|z\|$ 
and $\tQ\tz =x$.

Next let $(\tz_i)$ be as in the statement of c). 
Since $\|\tQ\|=1$ we need only prove the right hand inequality in c). 
Let $(a_i) \in c_{00}$ and choose $k\le m$ so that 
$$\trivert{\sum a_i \tz_i} = \|\tQ P_{[k,m]}^{\tbE} (\sum a_i \tz_i)\|$$

For all $i$, $P_{[k,m]}^{\tbE} \tz_i$ is $\tz_i$ or 0 except possibly for 
2 values of $i$, denoted by $i_0 \le i_1$. 
Thus 
\begin{equation*}
\|\tQ P_{[k,m]}^{\tbE} (\sum a_i \tz_i)\|
 \le |a_{i_0}| \, \trivert{\tz_{i_0}}
+ \|\sum_{i\in (i_0,i_1)} a_i \tQ (\tz_i)\| 
+ |a_{i_1}|\, \trivert{\tz_{i_1}}
 \le \frac{3\oK}{a} \|\sum a_i \tQ (\tz_i)\|\ ,
\end{equation*}
using that $(\tQ (\tz_i))$ has projection constant $\oK$ and is bounded 
below in norm by $a$.

It remains only to prove that $(\tE_i)$ is shrinking. 
Let $(\tz_i)$ be a normalized block sequence of $(\tE_i)$ in $\tZ$. 
Then $(\tQ\tz_i)$ is a bounded sequence in $X$. 
Moreover since $c_{00} (\oplus_{i=1}^\infty E_i^*) \cap X^*$ is dense in $X^*$,
$(\tQ \tz_i)$ is pointwise null on $X^*$ and hence weakly null in $X$.
We pass to a subsequence which we relabel as $(\tQ \tz_i)$ which is either 
norm null or satisfies $\inf_i \|\tQ z_i\| >0$ and is basic. 
In the latter case $(\tz_i)$ is weakly null by part~c). 
In the former case, given $n$ we can find a subsequence $(\tz_{i_j})_{j=1}^n$ 
with $\|\tQ \tz_{i_j}\| <\frac1n$ if $1\le j\le n$. 
Then if $\tz = \frac1n \sum_{j=1}^n \tz_{i_j}$, for some $k\le m$ 
and $j_0 \le j_1$ 
\begin{equation*}
\trivert{\tz} 
  = \|\tQ P_{[k,m]}^{\tbE} \tz\|
  \le   \Big\| \frac1n \tQ P_{[k,m]}^{\tbE} \tz_{i_{j_0}} \Big\| 
+ \Big\|\frac1n \sum_{j\in (j_0,j_1)} \tQ \tz_{i_j} \Big\| 
+  \Big\|\frac1n \tQ P_{[k,m]}^{\tbE} \tz_{i_{j_1}} \Big\| < \frac3n\ .
\end{equation*}
Thus in any case every normalized 
block sequence $(\tz_i)$ admits a convex block 
basis which is norm null and hence $(\tz_i)$ is weakly null and so 
$(\tE_i)$ is shrinking.
\end{proof}

We shall produce $A<\infty$ and a blocking $\tbH=(\tH_n)$ of $(\tE_n)$ with the 
following property. 
Let $x\in S_X$. 
Then there exists $\tz = \sum \tz_n \in \tZ$, with $\tz_n \in \tH_n$ 
for all $n$, so that if $(\tw_n)$ is any blocking of $(\tz_n)$ then 
$(\sum \trivert{\tw_n}^p)^{1/p} \le A$.
Moreover $\|\tQ \tz -x\| <\frac12$. 
Thus if $\tZ_p = \tZ_p (\tbH)$ then $\tQ :\tZ_p \to X$ remains an onto map. 
Moreover $\tZ_p$ is reflexive by 
Proposition~\ref{shrinkingFDD} and Lemma~\ref{part-a)} and $(\tH_n)$ is an FDD
for $\tZ_p$ satisfying 1-$(p,1)$-estimates. 
To accomplish this we need 

\begin{lem}\label{assume (2) holds}
Assume that \eqref{eq2} holds for our original map $Q:Z\to X$. 
Then there exist integers $0= N_0 < N_1 <\cdots$ so that if we define 
blockings $C_n = [E_i]_{i\in (N_{n-1},N_n]}$ and 
$D_n = [F_i]_{i\in (N_{n-1},N_n]}$ we have the following. 
Set for $n\in\nat$, 
\begin{equation*}
\begin{split}
& R_n = \left\{ i\in \nat :  N_n \ge i > \frac{N_n + N_{n-1}}2\right\}\ ,\\
& L_n = \left\{ i\in \nat : N_{n-1} < i\le  \frac{N_n + N_{n-1}}2\right\}\ ,\\
& C_{n,R} = [E_i]_{i\in R_n}\ \text{ and }\ 
C_{n,L} = [E_i]_{i\in L_n}\ .
\end{split}
\end{equation*}
Let $x\in S_X$, $i<j$, $\ep>0$ and assume that 
$\| P_{[1,i] \cup [j,\infty)}^\bD x \| <\ep$. 
Then there exists $z\in B_Z$ with $z\in [C_{i,R} \cup (C_\ell)_{\ell\in (i,j)}
\cup C_{j,L}]$ and 
$\|Qz - x\| < K[2\ep +\delta_i]$.
\end{lem}

\begin{proof} 
By \cite{Jo}  (see Lemma~4.3a \cite{OS})
we can choose  $0= N_0 <N_1 <\cdots$ so that if $z\in B_Z$ 
with $z= \sum z_j$, $z_j \in E_j$ for all $j$, then for $n\in\nat$ there 
exist $r_n\in R_n$ and $\ell_n \in L_n$ with $\|z_{r_n}\| <\ep_n$ and 
$\|z_{\ell_n} \| <\ep_n$. 
Define $C_n$ and $D_n$ as in the statement of the lemma and let $x\in S_X$ 
and $i<j$ with $\|P_{[1,i]\cup [j,\infty)}^\bD x\| <\ep$. 
Let $\|\bar z\| =1$ with $Q\bar z =x$ and $\bar z = \sum z_j$, 
$z_j\in E_j$ for all $j$. 
Choose $r_i\in R_i$ and $\ell_j \in L_j$ with $\|z_{r_i}\| <\ep_i$ and 
$\|z_{\ell_j}\| <\ep_j$. 
Let $z= \sum_{s\in (r_i,\ell_j)} z_s$. 
Thus $\|z\| \le 1$ and 
$z\in [C_{i,R} \cup (C_\ell)_{\ell \in (i,j)} \cup C_{j,L}]$.

Now 
\begin{itemize}
\item[{i)}] $\|P_{[1,r_i)\cup [\ell_j,\infty)}^\bF Qz\| 
< \ep_{r_i}  + \ep_{\ell_j -1}$ by \eqref{eq2}. 
\end{itemize}
Also if $w = \bar z -z = \sum_{s\notin (r_i,\ell_j)} z_j$ then 
again by \eqref{eq2} and our choice of $r_i$ and $\ell_j$, 
\begin{itemize}
\item[{ii)}] $\|P_{[r_i,\ell_j)}^\bF Qw\| 
= \| P_{[r_i,\ell_j)}^\bF Q 
(\sum_{s<r_i} z_s + z_{r_i} + z_{\ell_j} + \sum_{s>\ell_j} z_s)\| 
< K[\ep_{r_i} + \ep_i + \ep_j + \ep_{\ell_j +1}]$. 
\end{itemize}
{From} our hypothesis on $x$,
\begin{itemize}
\item[{iii)}] $\|P_{[1,r_i) \cup [\ell_j,\infty)}^\bF x\| < 2K\ep$.
\end{itemize}
Combining  i)--iii) we have, since $Qw = x-Qz$, 
\begin{equation*}
\begin{split}
\|Qz -x\| 
& \le \|P_{[1,r_i) \cup [\ell_j,\infty)}^\bF (Qz-x)\| 
    + \|P_{[r_i,\ell_j)}^\bF (Qz-x)\| \\
& < \ep_{r_i} + \ep_{\ell_j-1} + 2K\ep + K [\ep_{r_i} +\ep_i +\ep_j 
+ \ep_{\ell_j +1}]\\
& < K [2\ep + 6\ep_i] < K [2\ep +\delta_i]
\end{split}
\end{equation*}
since by \eqref{eq1} $6\ep_i <\delta_i$.
\end{proof}

We let $(C_n)$ and $(D_n)$ be the blockings given by 
Lemma~\ref{assume (2) holds}. 
Finally we block again using Proposition~\ref{OSprop4.4} for $(\delta_i)$ 
and $(D_n)$ to obtain $G_n = [D_i]_{i\in (k_{n-1},k_n]}$ for some 
$0= k_0 < k_1 <\cdots$. 
We set for $n\in\nat$, $H_n = [C_i]_{i\in (k_{n-1},k_n]}$. 

Let $x\in S_X$. 
Then by Proposition~\ref{OSprop4.4} there exists $(x_i)\subseteq X$ 
with $x= \sum x_i$ and for all $i\in \nat$ there exists $t_i \in (k_{i-1},k_i)$
so that 
\begin{itemize}
\item[{a)}] either $\|x_i\| <\delta_i$ or 
$\| P_{(t_{i-1},t_i)}^\bD x_i-x_i\| < \delta_i \|x_i\|$
\item[{b)}] $\|x_i\| < K+1$.
\end{itemize}

Let $B = \{i\in\nat : \|P_{(t_{i-1},t_i)}^\bD x_i-x_i \|< \delta_i \|x_i\|\}$ 
and $y =  \sum_{i\in B} x_i$. 
Then $\|x-y\| \le \sum_{i\notin B} \|x_i\| < \sum \delta_i < 1/4$ by 
\eqref{eq1}. 
For $i\in B$ set $\bar x_i = x_i/\|x_i\|$. 
{From} Lemma~\ref{assume (2) holds} there is a block sequence $(z_i)_{i\in B}$ 
of $(E_n)$ in $B_Z$ with 
\begin{equation}\label{eq3} 
\|Qz_i - \bar x_i\| < K[2\delta_i + \delta_{t_{i-1}}] < 3K\delta_i\ .
\end{equation}
Indeed the lemma yields that 
\begin{equation*}
z_i \in [C_{t_{i-1},R} \cup (C_\ell)_{\ell\in (t_{i-1},t_i)} \cup 
C_{t_i,L} ]
\end{equation*}
which ensures that the $z_i$'s are a block sequence.
{From} our choice of $(\delta_i)$ (right before \eqref{eq1}) and \eqref{eq3} we have that 
$(\tQ \tz_i)_{i\in B}$ is $2K$-basic, is 2-equivalent to $(\bar x_i)_{i\in B}$,
and $(\bar x_i)_{i\in B}$ satisfies $2C$-$(p,1)$-estimates.

{From} Proposition~\ref{shrinkingFDD} c) we have that 
\begin{equation}\label{eq4}
\begin{split}
\Big\|\sum_{i\in B} a_i \tQ \tz_i \Big\|  &\le \Bigtrivert{\sum_{i\in B} a_i \tz_i}
 \le \frac{3(2K)}{\inf\limits_{j\in B} \|\tQ\tz_j\|} 
\Big\|\sum_{i\in B} a_i \tQ \tz_i\Big\| \\
&\le 7K \Big\|\sum_{i\in B} a_i \tQ \tz_i\Big\| 
 \le 14 K \Big\|\sum_{i\in B} a_i \bar x_i\Big\|\ .
\end{split}
\end{equation}
(We have used that $\inf_{j\in B} \|\tQ\tz_j\| > 1-3K\delta_1$ from \eqref{eq3} 
and $\frac6{(1-3K\delta_1)} < 7$ from \eqref{eq1}.)

Let $\tz = \sum_{i\in B} \|x_i\| \tz_i$. 
Then from \eqref{eq4}, $\tz \in \tZ$ and moreover since $y = \sum_{i\in B} 
\|x_i\| \bar x_i$, 
$$\|\tQ\tz -y\| 
\le \sum_{i\in B} \|x_i\|\, \|\tQ \tz_i - \bar x_i\| 
< \sum_{i\in B} (K+1) 3K\delta_i < \frac14\ \text{ (by \eqref{eq1})}\ .$$
Thus $\|\tQ \tz-x\|  < 1/2$.

Finally we show that $\trivert{\tz}_{_{Z_p}} \le A$ where $A= 70CK$ and 
$\trivert{\cdot}_{_{Z_p}}$ denotes the norm of $\tZ_p (\tbH)$.

Write $\tz = \sum \tw_i$ where $(\tw_i)$ is any blocking of 
$(\|x_i\| \tz_i)_{i\in B}$.
Say $\tw_j = \sum_{i\in I_j} \|x_i\|\tz_i$ where $ I_1 <I_2<\cdots$ is any 
partition of $B$. 
Then by \eqref{eq4} if $y_j = \sum_{i\in I_j} x_i$, 
\begin{equation}\label{eq5} (\sum \trivert{\tw_j}^p)^{1/p} 
\le 14 K (\sum \|y_j\|^p )^{1/p} 
< 14 K \cdot 2C\|y\| < 35 CK
\end{equation}
since $(\bar x_i)$ satisfies $2C$-$(p,1)$-estimates and 
$\|y\| < 5/4$.

It remains to show that if we write $\tz = \sum \th_n$ where $\th_n\in\tH_n$
for all $n$ and $(\tg_n)$ is any blocking of $(\th_n)$ then 
\begin{equation}\label{eq6}
(\sum \trivert{\tg_n}^p)^{1/p} \le 70 CK \equiv A\ .
\end{equation}
As in the proof of a) there exists a blocking $(\tw_i)$ of 
$(\|x_i\|\tz_i)_{i\in B}$ with 
$\tg_n = P_{(j_{n-1},j_n]}^\bH (\tw_n + \tw_{n+1})$ for some $0=j_0 <j_1<\cdots$
and so $\trivert{\tg_n} \le \trivert{\tw_n} + \trivert{\tw_{n+1}}$.
Thus \eqref{eq6} follows.

This completes the proof of Theorem~\ref{main-central}.\qed
\medskip

We need some last preliminary results before proving Theorem~\ref{mainthm}. 

\begin{lem}\label{preliminary results}
Let $X$ be a reflexive Banach space and let $1\le q<\infty$ and 
$\frac1q + \frac1{q'} =1$. 
If $X$ satisfies $(\infty,q)$-tree estimates then $X^*$ 
satisfies $(q',1)$-tree estimates.
\end{lem}

\begin{proof} The case $q=1$ is trivial so we assume $q>1$.
By \cite{Z} we may assume that $X\subseteq Z$, a reflexive space with a 
bimonotone $\FDD\ (E_n)$. 
Note that if $(f_n) \subseteq S_{X^*}$ is weakly null then there exists 
a subsequence $(f_{n_i})$ of $(f_n)$ and a weakly null sequence 
$(x_i)\subseteq S_X$ with $\lim_i f_{n_i} (x_i)=1$. 
Indeed let $(y_n) \subseteq S_X$ with $f_n(y_n)=1$ for all $n$. 
Choose a subsequence $(y_{n_i})$ which converges weakly to some $y\in X$. 
Then $f_{n_i} (y_{n_i} -y) \to1$. 
Since $(E_n)$ is bimonotone, $\|y_{n_i}-y\|\to 1$ and thus we may take 
$x_i = (y_{n_i} -y)/\|y_{n_i} -y\|$.

Now let $(f_\alpha)_{\alpha \in T_\infty}$ be a weakly null tree in $S_{X^*}$. 
Using the above observation by successively replacing the successors of a 
given node by a subsequence we obtain a full subtree 
$(g_\alpha)_{\alpha \in T_\infty}$ of $(f_\alpha)_{\alpha\in T_\infty}$ and 
a weakly null tree $(x_\alpha)_{\alpha\in T_\infty}$ in $S_X$  so that for all 
$\alpha \in T_\infty \cup \{\emptyset\}$, $f_{\alpha,i} (x_{\alpha,i})\to1$ 
as $i\to\infty$.

Let $\ep_i\downarrow 0$. 
By again replacing each successor sequence of nodes by a subsequence we 
obtain two full subtrees $(f'_\alpha)_{\alpha \in T_\infty}$ and 
$(x'_\alpha)_{\alpha \in T_\infty}$ satisfying: 
For all branches $(\alpha_i)_{i=1}^\infty$ of $T_\infty$ and for all 
$i,j\in\nat$ with $i\ne j$, 
$$|f'_{\alpha_i} (x'_{\alpha_j}) | < \ep_{\max (i,j)}
\quad\text{and}\quad
f'_{\alpha_i} (x'_{\alpha_i}) > \frac12\ .$$
Let  $(x'_{\alpha_i})_{i=1}^\infty$ be a branch of 
$(x_\alpha)_{\alpha\in T_\infty}$ satisfying 
$C$-$(\infty,q)$-estimates. 
Let $(b_i)_{i=1}^\infty \in S_{\ell_{q'}}$ and choose $(a_i)_{i=1}^\infty 
\in S_{\ell_q}$ with $1= \sum_{i=1}^\infty a_i b_i$. 
Then $\|\sum a_j x'_{\alpha_j}\| \le C$  and so 
\begin{equation*}
\begin{split}
C\Big\| \sum b_i f'_{\alpha_i}\Big\| 
& \ge (\sum b_i f'_{\alpha_i})(\sum a_j x'_{\alpha_j} )\\
& \ge \sum_{i=1}^\infty a_i b_i f'_{\alpha_i} (x'_{\alpha_i})
- \sum_{i=1}^\infty \sum_{j\ne i} |f'_{\alpha_i} (x'_{\alpha_j})|\\
& > \frac12 - \sum_{i=1}^\infty [i\ep_i + \sum_{j>i} \ep_j] > \frac14
\end{split}
\end{equation*}
provided that the $\ep_i$'s were taken sufficiently small.
\end{proof}

\begin{prop}\label{bimonotone FDD}
Let $\bF= (F_i)$ be a bimonotone $\FDD$ for a Banach space $Z$ and let 
$1<q\le p<\infty$ and $C<\infty$. 
If $(F_i)$ satisfies $C$-$(\infty,q)$-estimates in $Z$ then $(F_i)$ is 
a bimonotone $\FDD$ for $Z_p (\bF)$ satisfying $C$-$(p,q)$-estimates in 
$Z_p(\bF)$.
\end{prop}

\begin{rem}
Prus (\cite[lemma 3.5]{P} obtained this result with weaker estimates. 
As written this result is stated with $C=1$ in \cite{JLPS}.
The clever proof we present was shown to us by W.B.~Johnson and G.~Schechtman.
\end{rem}

\begin{proof} 
Let $z= \sum z_i \in c_{00} (\oplus_{i=1}^\infty F_i)$ with $z_i\in F_i$ 
for all $i$. 
Let $k\in \nat$ and $0= n_0 <n_1 <\cdots < n_k \le\infty$. 
For some choice of $\ell$ and $0 = m_0 < m_1 <\cdots < m_\ell$ we have 
\begin{equation*}
\begin{split}
\Big\|\sum z_i\Big\|_{_{Z_p}}
& = \Big( \sum_{j=1}^\ell \Big\|\sum_{i=m_{j-1}+1}^{m_j} z_i\Big\|^p\Big)^{1/p}\\
& = \bigg[ \sum_{j=1}^\ell \Big\|\sum_{s=1}^k P_{(n_{s-1},n_s]}^\bF 
\Big( \sum_{i=m_{j-1}+1}^{m_j} z_i\Big) \Big\|^p \bigg]^{1/p}\\
& \le C \bigg[ \sum_{j=1}^\ell \bigg( \sum_{s=1}^k \|P_{(n_{s-1},n_s]}^\bF 
\Big( \sum_{i=m_{j-1}+1}^{m_j} z_i\Big)\|^q\bigg)^{p/q}\bigg]^{1/p}\\
&\le C\bigg[ \sum_{s=1}^k \bigg( \sum_{j=1}^\ell \| P_{(n_{s-1},n_s]}^\bF 
\Big( \sum_{i=m_{j-1} +1}^{m_j} z_i\Big) \|^p\bigg)^{q/p}\bigg]^{1/q}\ ,\\
&\qquad\qquad \text{by the ``reverse triangle inequality'' in } \ \ell_{p/q}^k\\
& = C\bigg[ \sum_{s=1}^k \bigg( \sum_{j=1}^\ell \|P_{(m_{j-1},m_j]}^\bF 
(P_{(n_{s-1},n_s]}^\bF z) \|^p \bigg)^{q/p} \bigg]^{1/q}\\
& \le C\bigg( \sum_{s=1}^k \|P_{(n_{s-1},n_s]}^\bF x\|_{_{Z_p}}^q\bigg)^{1/q}\ .
\end{split}
\end{equation*}
\end{proof}

For the next lemma and the proof of Theorem~\ref{mainthm} we adopt 
the convention that for $1\le p\le \infty$, $p'$ is defined by 
$1/p + 1/p' =1$.

\begin{lem}\label{Prus}
\cite{P}
Let $(E_i)$ be an $\FDD$ for a reflexive space $Z$. 
Let $1\le q\le p\le\infty$. 
Then $(E_i)$ satisfies $(p,q)$ estimates iff $(E_i^*)$ satisfies 
$(q',p')$ estimates.
\end{lem}

\begin{proof}[Proof of Theorem~\ref{mainthm}]
$\quad$

a) $\Rightarrow$ b).
Let $X$ be reflexive and satisfy $(p,q)$-tree estimates. 
By Lemma~\ref{preliminary results} $X^*$ satisfies $(q',1)$-tree estimates
(note that the case  $q=\infty$ is not possible, since $X$ is reflexive). 
Thus by Theorem~\ref{main-central}, $X^*$ is a quotient of a reflexive space 
$Z^*$ with an FDD $\bF^* = (F_n^*)$ satisfying 1-$(q',1)$-estimates. 
Hence by Theorem~\ref{main-central} $X$ embeds into $Z_p(\bG)$
for some blocking $\bG$ of $\bF$ and, 
by Lemma~\ref{Prus} and  
Proposition~\ref{bimonotone FDD}, $(G_n)$ satisfies $(p,q)$-estimates
 (again the case $p=1$ is not possible).

b) $\Rightarrow$ c) If $X$ satisfies b) then $p>1$ and, since the case $p=\infty$
 is trivial, by Theorem~\ref{main-central}, then $X$ is a quotient 
of a reflexive space with an FDD satisfying $(p,1)$-estimates. 
Thus by Lemma~\ref{Prus} $X^*$ is a subspace of a reflexive space with an 
FDD satisfying $(\infty, p')$-estimates. 
By Lemma~\ref{preliminary results} $X^*$ satisfies $(q',1)$-tree estimates 
and thus, by a) $\Rightarrow$ b), 
$X^*$ embeds into a reflexive space $Z^*$ with an FDD satisfying 
$(q',p')$-estimates. 
Hence $X$ is a quotient of $Z$, a reflexive space with an FDD satisfying 
$(p,q)$-estimates, again using Lemma~\ref{Prus}.

c) $\Rightarrow$ a). 
We assume that $X$ is a quotient of a reflexive space $Z$ having an FDD 
satisfying $(p,q)$-estimates. 
Thus $X^* \subseteq Z^*$ which by Lemma~\ref{Prus} has an FDD satisfying 
$(q',p')$-estimates. By b) $\Rightarrow$ c) we have that $X^*$ is a quotient
 of a reflexive space having an FDD satisfying $(q',p')$-estimates.
 Thus c) follows by duality using   Lemma~\ref{Prus}.
\end{proof}

\begin{remark}\label{equivalences}
The following equivalences can be added to Theorem~\ref{mainthm}. 
\begin{itemize}
\item[{d)}] $X$ is isomorphic to a subspace of a quotient of a reflexive 
space $Z$ having an FDD which satisfies $(p,q)$-estimates.
\item[{e)}] $X^*$ satisfies $(q',p')$-tree estimates.
\item[{f)}] $X$ is isomorphic to a subspace of a reflexive space $Z$ 
having an FDD which satisfies 1-$(p,q)$-estimates.
\end{itemize}
Indeed f) follows from the proof of a) $\Rightarrow$ b) in 
Theorem~\ref{mainthm} since $Z$ has an $\FDD\ (F_i)$ satisfying 
1-$(\infty,q)$-estimates and so by Lemma~\ref{Prus}, 
$Z_p(\bF)$ satisfies 1-$(p,q)$-estimates. 
Thus we obtain a solution to a question raised in \cite{JLPS} after the 
statement of Proposition~\ref{preliminary results}. 
We refer the reader to \cite{JLPS} for the relevant definitions.
\end{remark}

\begin{cor}\label{X reflexive}
Let $X$ be a reflexive Banach space and $1< q\le p <\infty$.
The following are equivalent
\begin{itemize}
\item[{a)}] $X$ embeds into a reflexive space $Z$ having an $\FDD$ satisfying 
$1$-$(p,q)$-estimates. 
\item[{b)}] $X$ can be renormed to be asymptotically uniformly smooth of 
power type $q$ and $X$ can be renormed to be asymptotically uniformly 
convex of power type $p$.
\item[{c)}] $X$ can be renormed so as to be both asymptotically uniformly 
smooth of power type $q$ and asymptotically uniformly convex of power 
type $p$.
\end{itemize}
\end{cor}

The following remains open.

\begin{prob}\label{problem}
Let $X$ be a uniformly convex separable Banach space. 
Does there exist a uniformly convex space $Z$ with an FDD (or even a basis) 
so that $X$ embeds into $Z$?
\end{prob}


\begin{thebibliography}{MM}
\frenchspacing

\bibitem[B]{B}  
J. Bourgain,
{\em On separable Banach spaces, universal for all separable reflexive spaces},
Proc. Amer. Math. Soc. {\bf79} (1980), no.2, 241--246.

\bibitem[E]{E} 
P. Enflo, 
{\em Banach spaces which can be given an equivalent uniformly convex norm},
Proceedings of the International Symposium on Partial Differential Equations 
and the Geometry of Normed Linear Spaces, Jerusalem, 1972, 
Israel J. Math. {\bf13} (1972), 281--288 (1973).

\bibitem[GG]{GG} 
V.I. Gurari\u i and N.I. Gurari\u i, 
{\em Bases in uniformly convex and uniformly smooth Banach spaces}, 
Izv. Akad. Nauk SSSR Ser. Mat. {\bf35} (1971), 210--215 (Russian).

\bibitem[J]{J} 
Robert C. James, 
{\em Super-reflexive spaces with bases}, 
Pacific J. Math. {\bf 41} (1972), 409--419.


\bibitem[Jo]{Jo} W.B. Johnson,
{\em On quotients of $L_p$ which are quotients of $l_p$}, 
Compositio Math. {\bf34} (1977), no.1, 69--89.

\bibitem[JZ1]{JZ1} W.B. Johnson and M. Zippin,
{\em On subspaces of quotients of $(\sum G_n)_{l_p}$ and $(\sum G_n)_{c_0}$}, 
Proceedings of the International Symposium on Partial Differential Equations 
and the Geometry of Normed Linear Spaces, Jerusalem, 1972, 
Israel J. Math. {\bf13} (1972), 311--316 (1973).

\bibitem[JZ2]{JZ2} W.B. Johnson and M. Zippin,
{\em Subspaces and quotient spaces of $(\sum G_n)_{l_p}$ and 
$(\sum G_n)_{c_0}$},
Israel J. Math. {\bf17} (1974), 50--55.

\bibitem[JLPS]{JLPS} 
W.B. Johnson, J. Lindenstrauss, D. Preiss, and G. Schechtman,
{\em Almost Fr\'echet differentiability of Lipschitz mappings between 
infinite-dimensional Banach spaces}, 
Proc. London Math. Soc. (3) {\bf 84} (2002), no.3, 711--746.

\bibitem[OS]{OS} 
E. Odell and Th. Schlumprecht,
{\em Trees and branches in Banach spaces},  
Trans. Amer. Math. Soc. {\bf354} (2002), no.10, 4085--4108.

\bibitem[Pi]{Pi} 
G. Pisier, 
{\em Martingales with values in uniformly convex spaces}, 
Israel J. Math. {\bf20} (1975), no.3--4, 326--350.

\bibitem[P]{P} S. Prus,
{\em Finite-dimensional decompositions with $p$-estimates and universal 
Banach spaces} (Russian summary), 
Bull. Polish Acad. Sci. Math. {\bf31} (1983), no.5--8, 281--288 (1984).

\bibitem[Z]{Z} M. Zippin, 
{\em Banach spaces with separable duals}, 
Trans. Amer. Math. Soc. {\bf310} (1988), no.1, 371--379.

\end{thebibliography}
\end{document}